


\input amstex
\documentstyle{amsppt}
\hsize28cc\vsize42cc 
\hoffset18mm\voffset8mm
\pageno 1 %
\TagsOnLeft
\loadbold

\topmatter

\vbox{\noindent\sevenrm\baselineskip9pt 
PUBLICATIONS DE L'INSTITUT MATH\'EMATIQUE\newline
\hglue4em Nouvelle s\'erie, tome 67(81) (2000), 14--30}\vglue20mm

\title               UN TH\'EOR\`EME DE MOYENNE
                     POUR LES SOMMES D'EXPONENTIELLES.
                     APPLICATION \`A L'IN\'EGALIT\'E DE WEYL
\endtitle

\author  O. Robert et P. Sargos
\bigskip\endauthor

 \leftheadtext{Robert et Sargos}
\rightheadtext{Un th\'eor\`eme de moyenne pour les sommes d'exponentielles}

\address\eightrm   Institut Elie Cartan   \hfill (Re\c cu le  07 10 1998)
  \break\indent    Universit\'e Henri Poincar\'e - Nancy I
\newline\indent    BP 239
\newline\indent    54 506 Vandoeuvre-l\`es-Nancy Cedex
\smallskip
{\eighttt          Robert\@iecn.u-nancy.fr
\newline\indent    Sargos\@iecn.u-nancy.fr}
\endaddress

\dedicatory        Communicated by Aleksandar Ivi\'c
\enddedicatory

\subjclass       Primary 11 L 07; Secondary 11 L 15 
\endsubjclass

\abstract
Nous \'etudions la moyenne de puissances sixi\`emes de certaines
sommes d'exponentielles selon une m\'ethode due \`a Bombieri et
Iwaniec [1]. Le r\'esultat s'applique \`a l'in\'egalit\'e de Weyl en
suivant une id\'ee due \`a Heath-Brown [3].
\endabstract

\endtopmatter

\def\d{d}
\def\e{e}
\def\ep{\varepsilon}
\def\bs{\boldsymbol}
\def\bk{\boldkey}
\def\r{\Bbb R}
\def\B{\Cal B}
\redefine\le{\leqslant}
\redefine\ge{\geqslant}

\vglue-2em

\document

\head{1. Introduction}\endhead
Soit $p$ un entier fix\'e, $p=3,4$ ou $5$. Nous cherchons \`a majorer la 
moyenne de puissances $2p$-i\`emes de sommes d'exponentielles:
$$
I_{2p}=\int_0^1\int_0^\lambda 
\bigg|\sum_{n=N+1}^{2N}\e(\alpha n^2+\gamma n^4)\bigg|^{2p}\d\alpha\,\d\gamma,
\leqno (1.1)
$$ 
avec $N$ entier $\ge 2$, $\lambda$ r\'eel positif et avec la notation 
$\e(x)=\e^{2i\pi x}$.

Un raisonnement heuristique facile montre que le r\'esultat attendu est:
$$
 I_{2p}\ll_\ep \lambda N^{p+\ep}+N^{2p-6+\ep}, \text{ \ pour tout \ } \ep>0 
\leqno(1.2) 
$$

Notre r\'esultat principal consiste \`a d\'emontrer (1.2) dans le cas
$p=3$.  Pour cela, nous suivons un sc\'enario d\^u \`a Bombieri et
Iwaniec [1]; nous nous d\'emarquons de [1] par de nombreuses
simplifications de d\'etail rendues possibles par le cadre plus simple
dans lequel nous nous pla\c cons.  Notre motivation principale pour ce
r\'esultat est un {\it{crit\`ere de la d\'eriv\'ee cinqui\`eme}} pour
les sommes d'exponentielles ; cette application est d\'etaill\'ee dans
[6].
Dans les cas $p=4$ et $p=5$, nous ne parvenons pas au r\'esultat
attendu (1.2), et nous nous contentons de d\'eduire du cas $p=3$ les
majorations suivantes:
$$
I_8 \ll_\ep \lambda N^{9/2+\ep}+N^{2+\ep},
\text{ \ pour tout \ } \ep>0 
$$ 
(pour fixer le id\'ees, signalons qu'un raisonnement direct
\'equivalent \`a la d\'emon\-stra\-tion du lemme de Hua, cf\. [7], aboutirait
seulement \`a la majoration: $I_8\ll_\ep\lambda N^{5+\ep}+N^{2+\ep}$)
et
$$
I_{10}\ll_\ep \lambda N^{49/8+\ep}+N^{4+\ep}, 
\text{ \ pour tout \ } \ep>0 
\leqno (1.3)
$$

Ce dernier r\'esultat s'applique \`a l'in\'egalit\'e de Weyl en suivant
une id\'ee due \`a Heath-Brown [3]. Plus pr\'ecis\'ement, soit $k$ un
entier fix\'e. Etant donn\'e le r\'eel $\alpha$ et l'entier $N$, on
pose:
$$
S_k (\alpha)=\sum_{n=1}^N \e(\alpha n^k).
$$

Supposons que $\alpha$ admette une approximation rationnelle
``g\'en\'erique" (en un sens \`a pr\'eciser dans chaque cas). Alors
l'in\'egalit\'e classique de Weyl s'\'ecrit, en posant $K=2^k$:
$S_k(\alpha)\ll_\ep N^{1-2/K+\ep}$, pour chaque $k\ge 2$ (cf. [4] ou
[7]). Dans le th\'eor\`eme 1 de [3], l'exposant $1-2/K$ devient
$1-8/3K$, pour $k\ge 6$, avec des conditions plus restrictives sur
$\alpha$.

Comme cons\'equence de (1.3), nous obtenons l'exposant $1-3/K$, pour
$k\ge 8$, avec de nouvelles restrictions sur $\alpha$.

\head{2. Notations}\endhead

L'\'ecriture $u\ll v $ signifie qu'il existe une constante absolue
$C>0$ telle qu'on ait $|u|\le C v$. l'\'ecriture $u\ll_k v $ signifie
que $C$ peut d\'ependre de $k$. L'\'ecriture $u\ll_\ep v $ sous-entend
que la majoration a lieu pour tout $\ep >0$. Si $u$ et $v$ sont
positifs, l'\'ecriture $u\asymp v$ signifie qu'on a \`a la fois 
$u\ll v$ et $v\ll u$.

La fonction caract\'eristique de l'intervalle $[a,b]$ est not\'ee
${\bk 1}_{[a,b]}$.

Enfin, le symbole $\square$  se place \`a la fin d'une
d\'emonstration pour signaler que celle-ci est termin\'ee, ou \`a la
fin d'un \'enonc\'e pour signaler que la d\'emonstration a \'et\'e
omise.

\head{3. Enonc\'e du r\'esultat principal}\endhead

Afin de formuler nos r\'esultats dans un cadre ad\'equat, nous
introduisons quelques notations.

Soit un r\'eel $N\ge 2$. Pour toute suite $(b_n)_{N<n \le 2N}$ de
nombres complexes, on pose:
$$
 \bigg|\sum_{n \sim N}b_n \bigg| =\max_{N<N_1 \le 2N} 
\bigg| \sum_{N<n\le N_1}b_n \bigg|. 
\leqno (3.1)
$$
De m\^eme, on d\'esigne par $\sum_{n \sim N}b_n $ la somme
$\sum_{N<n\le N_1}b_n $ pour laquelle le maximum ci-dessus est
r\'ealis\'e, $N_1$ \'etant choisi minimum en cas d'ambigu\"{\i}t\'e.

Dans toute la suite, $(a_n)_{N<n \le 2N}$ d\'esigne une suite de
nombres complexes de modules au plus \'egaux \`a $1$. Soit $\lambda>0$.
On pose enfin:
$$
 I_{2p}=\int_0^1\int_0^\lambda 
\bigg|\sum_{n\sim N}a_n \e(\alpha n^2+\gamma n^4)\bigg|^{2p}\d\alpha\,\d\gamma,
$$

Notre r\'esultat principal s'\'enonce ainsi:

\proclaim{Th\'eor\`eme 1} Avec les notations ci-dessus, on a la
majoration suivante:
$$
I_6\ll_\ep \lambda N^{3+\ep}+N^{\ep}.
$$
\endproclaim

Les deux r\'esultats suivants sont cons\'equence du th\'eor\`eme 1.

\proclaim{Th\'eor\`eme 2} Avec les notations ci-dessus, on a de m\^eme :
$$
I_8\ll_\ep \lambda N^{9/2+\ep}+N^{2+\ep} 
\leqno (3.2)
$$
\endproclaim

\proclaim{Th\'eor\`eme 3} Avec les notations ci-dessus, on a enfin :

$$
I_{10}\ll_\ep \lambda N^{49/8+\ep}+N^{4+\ep} 
$$
\endproclaim

\head
{4. Syst\`emes diophantiens et moyennes de sommes d'exponentielles}
\endhead

Dans cette section, nous \'etablissons un lemme qui, bien que
compl\`etement \'el\'ementaire, est essentiel dans la d\'emonstration
des th\'eor\`emes 1, 2 et 3. Il exprime que les bornes des int\'egrales
et du domaine de sommation dans (1.1) peuvent \^etre modifi\'ees \`a
volont\'e, moyennant un facteur multiplicatif. Le ph\'enom\`ene
s'explique par le lien qui existe entre $I_{2p}$ et le syst\`eme
diophantien:

$$
\aligned
 |n_1^2+\cdots+n_p^2-m_1^2-\cdots-m_p^2| &\le  \delta_1N^2 \\
 |n_1^4+\cdots+n_p^4-m_1^4-\cdots-m_p^4| &\le  \lambda_1N^4
\endaligned
\quad
 N<n_i,m_i\le 2N 
\leqno (4.1)
$$

avec $\delta_1$ et $\lambda_1$ r\'eels positifs ou nuls, dont le
nombre de solutions sera not\'e:
$$
{\Cal B}_{2p}(N,\delta_1,\lambda_1). 
\leqno (4.2)
$$
Ce lien est tout \`a fait classique (cf\. [4, chapitre 4]) et a \'et\'e
largement utilis\'e sous une forme semblable dans [1] et
[8]. Cependant, nous pr\'ef\'erons red\'emontrer compl\`etement le
r\'esultat pr\'ecis dont nous avons besoin:

\proclaim{Lemme 1}
Soit $p\ge 1$ et $N\ge 2$ deux entiers, avec $p$ fix\'e.
 Soient 
$c,d,\lambda,\mu,\Delta, \delta $ des r\'eels tels que $0<\Delta \le
\delta \le 1$, $0< \mu \le\lambda $.
Soient, pour $\alpha\in [c,c+\delta]$ et $\gamma\in [d,d+\lambda]$,
des fonctions $\varphi_{\alpha, \gamma}: 
[N,2N] \to \r $ de classe ${\Cal C}^1$ telles que 
$\varphi'_{\alpha, \gamma} \ll 1/N$.

Alors on a
$$
\align
\frac{1}{(\log N)^{2p}}\frac{1}{\delta \lambda}
&\int_c^{c+\delta}\int_{d}^{d+ \lambda}\bigg|
\sum_{n\sim N}a_n e(\alpha n^2 + \gamma n^4+ \varphi_{\alpha, \gamma}(n))
\bigg|^{2p}\d\alpha\,\d\gamma \\
&\ll_p \B_{2p}\Big(N,\frac{1}{\delta N^2},\frac{1}{\lambda N^4}\Big)
\tag 4.3\\
&\ll_p\frac {1}{\Delta \mu}\int_{- \Delta /2}^{\Delta /2}
\int_{-\mu /2}^{\mu /2} \bigg|\sum_{N<n\le 2N} e(\alpha n^2 + \gamma n^4)
\bigg|^{2p} \d\alpha\,\d\gamma \cr
\endalign
$$
\endproclaim 
{\it Remarque:} Si la fonction 
$$
 (\alpha,\gamma)\to \bigg|\sum_{n \sim
N}a_n e(\alpha n^2 + \gamma n^4+ \varphi_{\alpha, \gamma}(n))
\bigg|^{2p}
$$
n'est pas mesurable au sens de Lebesgue, l'int\'egrale du premier
membre de (4.3) doit \^etre pris au sens d'une int\'egrale
sup\'erieure.

{\it D\'emonstration. }
1) Pour chaque $(\alpha,\gamma)$ fix\'e, la fonction
$n\to\varphi_{\alpha,\gamma}(n)$ a une variation totale $\ll 1/N$. Par
le choix des notations, la sommation d'Abel s'\'ecrit:
$$
\sum_{n \sim
N}a_n e(\alpha n^2 + \gamma n^4+ \varphi_{\alpha, \gamma}(n))
\ll \bigg|\sum_{n \sim
N}a_n e(\alpha n^2 + \gamma n^4)\bigg| 
$$
Par ailleurs, pour tout $N_1\in]N,2N]$, on a 
$$
\multline
\sum_{N<n\le N_1}a_ne(\alpha n^2 + \gamma n^4)\\
=\int_{-1/2}^{1/2}\biggl(\sum_{N<n\le 2N}a_ne(\alpha n^2 + \gamma
n^4-\theta n)\biggr) \biggl(\sum_{N<\nu\le N_1}e(\nu \theta)\biggr)\,\d\theta 
\endmultline
$$
Posons ${\Cal L}(\theta)=\min(N,1/|\theta|)$. Par
l'in\'egalit\'e de H\"older, on a:
$$
\multline
\bigg|\sum_{n \sim N}a_ne(\alpha n^2 + \gamma n^4)\bigg|^{2p}\\
\ll_p
(\log N)^{2p-1}\int_{-1/2}^{1/2}\L(\theta)\bigg|\sum_{N<n\le 2N}a_ne(\alpha 
n^2 + \gamma n^4-\theta n)\bigg|^{2p}\d \theta.
\endmultline
$$
En int\'egrant par rapport \`a $\alpha $ et $\gamma$, on obtient
finalement
$$
\multline
\frac {1}{(\log N)^{2p}}\frac{1}{\delta \lambda 
}\int_c^{c+\delta}\int_{d}^{d+ \lambda}
\bigg|\sum_{n \sim N}a_n e(\alpha n^2 + \gamma n^4+ \varphi_{\alpha, \gamma}(n))
\bigg|^{2p}\d\alpha\,\d\gamma\\
\ll_p \frac{1}{(\log N)}\int_{-1/2}^{1/2}{\Cal L}(\theta)I(\theta)\,\d\theta,
\endmultline
$$

avec
$$
I(\theta)=\frac{1}{\delta \lambda 
}\int_c^{c+\delta}\int_{d}^{d+ \lambda}
\bigg|\sum_{N<n\le 2N}a_n'(\theta) e(\alpha n^2 + \gamma n^4)
\bigg|^{2p}\d\alpha\,\d\gamma 
$$
et $a_n'(\theta)=a_n\e (-\theta n)$.

2) Pour $\theta $ fix\'e, on va \'etablir la majoration:
$$
I(\theta)\ll {\Cal B}_{2p}\Bigl(N,\frac{1}{\delta N^2},\frac{1}{\lambda
N^4}\Bigr), 
\leqno (4.4)
$$
ce qui suffit \`a prouver la premi\`ere in\'egalit\'e dans (4.3). 

Soit la fonction $\rho(y)=\Big(\frac{\sin (\pi y)}{\pi y}\Big)^2 $ (avec
la convention $\rho(0)=1$), dont la transform\'ee de Fourier vaut
$$
\hat{\rho}(t)=\int_{-\infty}^{+\infty}\rho(y)\e(-ty)\d y=(1-|t|)^{+},
$$
et qui v\'erifie: 
$$
{\bk 1}_{[-1/2,1/2]}(y)\ll \rho(y).
$$
 On revient
\`a $I(\theta)$. Par changement de variable, on \'ecrit:
$$
I(\theta)=\int_{-1/2}^{1/2}\int_{-1/2}^{1/2}
\bigg|\sum_{N<n\le 2N}a_n''(\theta) e(\alpha\delta n^2 + \gamma\lambda n^4)
\bigg|^{2p}\d\alpha\,\d\gamma, 
$$
avec $a_n''(\theta)=a_n\e((c+\delta/2)n^2+(d+\lambda/2)n^4-\theta n)$,
d'o\`u
$$
I(\theta)\ll\int_{\r}\int_{\r}\rho(\alpha)\rho(\gamma)
\bigg|\sum_{N<n\le 2N}a_n''(\theta) e(\alpha\delta n^2 + \gamma\lambda n^4)
\bigg|^{2p}\d\alpha\,\d\gamma, 
$$
On d\'eveloppe maintenant la puissance $2p-$i\`eme et on int\`egre
terme \`a terme. Tenant compte du fait que $|a_n|\le 1$, on arrive
sans mal \`a la formule:
$$
I(\theta)\ll \sum_{{\bs n}}\sum_{{\bs m}}
\hat{\rho}\big(\delta (s_2({\bs n})-s_2({\bs m}))\big)
\hat{\rho}\big(\lambda (s_4({\bs n})-s_4({\bs m}))\big)
\leqno (4.5)
$$
o\`u on a pos\'e:

${\bs n}=(n_1,\dots,n_p)$ et
${\bs m}=(m_1,\dots,m_p)\in ]N,2N]^p$,
$s_2({\bs n})=n_1^2+\cdots+n_p^2$ et
$s_4({\bs n})=n_1^4+\cdots+n_p^4$. 
Il est clair que (4.5) implique (4.4).

3) Nous allons \'etablir la deuxi\`eme in\'egalit\'e de (4.3). On
remarque d'abord l'in\'egalit\'e suivante qui, bien que triviale,
constitue le coeur de la d\'emonstration du lemme 1:
$$
{\Cal B}_{2p}\Bigl(N,\frac{1}{\delta N^2},\frac{1}{\lambda N^4}\Bigr)
\le {\Cal B}_{2p}\Bigl(N,\frac{1}{\Delta N^2},\frac{1}{\mu N^4}\Bigr)
$$
Pour finir la d\'emonstration du lemme 1, il suffit de prouver que:
$$
{\Cal B}_{2p}\left(N,\frac{1}{\Delta N^2},\frac{1}{\mu N^4}\right)
\ll \frac{1}{\Delta\mu
}\int_{-\Delta/2}^{\Delta/2}\int_{-\mu/2}^{\mu/2}
\bigg|\sum_{N<n\le 2N} e(\alpha n^2 + \gamma n^4)
\bigg|^{2p}\d\alpha\,\d\gamma, 
\leqno (4.6)
$$
et, pour cela, de reprendre les calculs de la deuxi\`eme partie en sens 
inverse. Plus pr\'ecis\'ement, on a:
$$
\multline
{\Cal B}_{2p}\Bigl(N,\frac{1}{\Delta N^2},\frac{1}{\mu N^4}\Bigr)\\
=\sum_{{\bs n}}\sum_{{\bs m}} {\bk 1}_{[-1/\Delta,1/\Delta]}\big(s_2({\bs n})
-s_2({\bs m})\big){\bk 1}_{[-1/\mu,1/\mu]}\big(s_4({\bs n})-s_4({\bs m})\big)
\endmultline
\tag 4.7
$$
(avec ${{\bs n}}$ et $ {{\bs m}}\in ]N,2N]^p$).

Maintenant, on \'ecrit que 
$
{\bk 1}_{[-1/\Delta,1/\Delta]}(y)\ll\rho({\Delta y}/2)
$
et que 
$
\rho(y)=\smash{\int\limits_{\r}}\hat{\rho}(x)\e(xy)\d x,
$
ce qui montre que le deuxi\`eme membre de (4.7) est 
$$
\align
&\ll \int_{\r}\int_{\r}
\hat{\rho}(\alpha)\hat{\rho}(\gamma)
\bigg\{
\sum_{{\bs n}}\sum_{{\bs m}}
\e\left(\frac{\alpha\Delta}{2}(s_2({\bs n})-s_2({\bs m}))
+\frac{\gamma \mu}{2}(s_4({\bs n})-s_4({\bs m}))\right)\!\!
\bigg\}\d\alpha\,\d\gamma 
\\
&=\int_{\r}\int_{\r}
\hat{\rho}(\alpha)\hat{\rho}(\gamma)
\bigg|\sum_{N<n\le 2N} e(\alpha n^2 + \gamma n^4)
\bigg|^{2p}\d\alpha\,\d\gamma.
\endalign
$$
On majore $\hat{\rho}(y)$ par ${\bk 1}_{[-1,1]}(y)$, puis on fait un
changement de variable \'evident et on obtient exactement (4.6).
\hfill\hfill\qed

\head{5. Puissances quatri\`emes}\endhead

Le but de cette section est de majorer l'int\'egrale
$$
\int_0^{N^{-1/2}}\int_{-N^{-3}}^{N^{-3}}
\bigg|\sum_{n \sim
N}\e(\alpha n^2 + \gamma n^4)
\bigg|^{6}\d\alpha\,\d\gamma,
$$
ce qui fournit l'analogue du corollaire de [1]. Auparavant, nous
devons \'etablir un lemme auxiliaire sur les puissances quatri\`emes:

\proclaim{Lemme 2} Soient un entier $N\ge 2$ et un r\'eel $\Delta\ge
N^{-1}$. Alors:
$$
\int_0^{\Delta}\int_{-1/N^3}^{1/N^3}
\bigg|\sum_{n\sim N}a_n\e(\alpha n^2+\gamma n^4)\bigg|^4 \d\alpha\,\d\gamma
\ll \frac{\Delta(\log N)^5}{N}. 
\leqno (5.1)
$$
\endproclaim

 Soit ${\Cal J}$ l'int\'egrale du premier membre de (5.1). Le
lemme 1 s'applique ici:
$$
{\Cal J}\ll (\log N)^4 {\Cal B}_4\Big(N,\frac{1}{\Delta N^2},
\frac{1}{N}\Big)\le (\log N)^4 {\Cal B}_4\Big(N,\frac{1}{N},
\frac{1}{N}\Big), 
$$
la derni\`ere in\'egalit\'e \'etant \'evidente du fait que $\Delta \ge
\frac{1}{N}$. On rappelle que $ {\Cal B}_4= {\Cal B}_4\Big(N,\frac{1}{N},
\frac{1}{N}\Big)$ est le nombre de solutions du syst\`eme:
$$
\aligned
 |n_1^2+n_2^2-m_1^2-m_2^2| &\le N \\
 |n_1^4+n_2^4-m_1^4-m_2^4| &\le  N^3
 \endaligned\qquad
 N<n_i,m_i\le 2N 
$$
qu'on majore par le nombre de solutions du syst\`eme:
$$
\aligned
&a)\quad n_1^2+n_2^2 = m_1^2+m_2^2 +O(N)\\
&b)\quad n_1^4+n_2^4 = m_1^4+m_2^4 +O(N^3)
\endaligned\qquad
 N<n_i,m_i\le 2N 
\tag 5.2
$$
Si on \'el\`eve (5.2.a) au carr\'e et qu'on lui retranche (5.2.b), on
obtient:
$$
(n_1n_2)^2=(m_1m_2)^2+O(N^3),
$$ 
ce qui implique imm\'ediatement
$$
n_1n_2=m_1m_2+O(N),
\leqno (5.3)
$$ 
Dans le syst\`eme constitu\'e par (5.2.a) et (5.3), on pose 
$$
a=n_1+n_2, \enskip b=m_1+m_2,\enskip c=n_1-n_2,\enskip d=m_1-m_2.
$$
 On a montr\'e que ${\Cal
B}_4$ est major\'e par le nombre de solutions d'un syst\`eme du type:
$$
\alignedat2
&a)\quad a = b +O(1)&  &a\asymp N,\ b\asymp N,\\
&b)\quad\  c^2 = d^2 +O(N)\qquad && c\ll N,\ d\ll N
\endalignedat
\tag 5.4
$$
Supposons que $c$ soit d'un ordre de grandeur fix\'e $U$, avec $1\le
U\le N$. Alors (5.4.b) s'\'ecrit: 
$$
c-d\ll N/U,
$$
qui poss\`ede $O(N)$ solutions. En sommant sur les $O(\log N)$ valeurs
possibles de $U$, et en ajoutant le cas $c=0$, on a montr\'e que le
nombre de solutions de (5.4) est $\ll N^2\log N$, et de m\^eme pour
${\Cal B}_4$. \hfill\hfill\qed

Nous sommes maintenant en mesure d'\'etablir le:

\proclaim{Lemme 3} Pour tout entier $N\ge 2$, on a:
$$
\int_0^{N^{-1/2}}\int_{-N^{-3}}^{N^{-3}}
\bigg|\sum_{n \sim
N}\e(\alpha n^2 + \gamma n^4)
\bigg|^{6}\d\alpha\,\d\gamma \ll (\log N)^6 
$$
\endproclaim

On d\'ecoupe le domaine d'int\'egration de la variable $\alpha $
en $O(\log N)$ intervalles:
$$
\int_0^{N^{-1/2}}\ \le \int_0^{16N^{-1}}
+\sum_{16N^{-1}\le \delta < N^{-1/2}}\int_{\delta}^{2\delta}\, 
\text{ \ avec \ } \delta = 2^kN^{-1},
$$
et on montre que chaque terme qui lui correspond est $\ll (\log N)^5$.
Par exemple, soit $\delta$ tel que $16N^{-1}\le \delta <N^{-1/2} ;$ on pose
$$
S(\alpha,\gamma)=\sum_{n\sim N}\e(\alpha n^2+\gamma n^4)
\quad
\text{et}
\quad
I_{\delta}=\int_{\delta}^{2\delta}\int_{-N^{-3}}^{N^{-3}}\big|S(\alpha
,\gamma)\big|^6\d\alpha\,\d\gamma.
$$
On \'ecrit alors:
$$
I_{\delta}\ll \Big(\max_{\delta\le \alpha\le 2\delta}
\big|S(\alpha,\gamma)\big|\Big)^2 \int_{0}^{2\delta}\int_{-N^{-3}}^{N^{-3}}\big|S(\alpha
,\gamma)\big|^4\d\alpha\,\d\gamma.
$$
Comme on a choisi $ \delta\ge 16N^{-1}$, on peut appliquer
l'in\'egalit\'e de Van der Corput (cf\. [2, Th\'eo\-r\`eme~2.2]) pour
obtenir $S(\alpha,\gamma)^2\ll N^2\delta$, ce qui montre, \`a l'aide
du Lemme 2, qu'on a:
$$
 I_{\delta}\ll N\delta^2(\log N)^5\le (\log N)^5. \hfill\hfill\qed
$$

\head{6. Un calcul sur la transformation B}\endhead

Comme dans la d\'emonstration de Bombieri et Iwaniec (cf\. [1, section
3]), nous avons besoin de calculs pr\'ecis sur la transformation B de
Van der Corput. Ces calculs ont, depuis, \'et\'e \'ecrits dans un cadre
g\'en\'eral (cf\. [5]). Nous nous contentons ici de rappeler le
r\'esultat dont nous avons besoin.
 Soit $f:[N,AN]\to \r $ une fonction de classe ${\Cal C}^4$ telle
que:
$$
\lambda_2\le f''(x)\ll\lambda_2,\enskip f'''(x)\ll\lambda_2/N,\enskip
f^{(4)}(x)\ll\lambda_2/N^2
\leqno (6.1) 
$$
pour $x\in[N,AN]$, o\`u $\lambda_2$ d\'esigne un r\'eel positif.

La transformation B de Van der Corput (cf\. [2, Lemme 3.6], ou
[4, chapitre 3, Th\'eor\`eme 10]) s'\'ecrit:
$$
\sum_{N<n\le 2N}\e(f(n))=\e^{i\pi /4}
\sum_{\alpha\le \nu\le\beta}\frac{\e(f^{*}(\nu))}{f''(z(\nu))}+
O(\log (2+N\lambda_2))+O(\lambda_2^{-1/2}), 
\leqno (6.2)
$$
avec les notations suivantes: $\alpha =f'(N), \beta = f'(AN)$,
$$
z(y)=\big(f'\big)^{-1}(y), \text{ \ pour \ } \alpha\le y \le \beta
\leqno (6.3)
$$
et enfin
$$
f^{*}(y)=f(z(y))-yz(y), \text{ \ pour \ } \alpha\le y \le \beta.
\leqno (6.4)
$$
Cette derni\`ere relation d\'efinit la fonction $f^{*}$ \`a partir de
$f$.

Soit maintenant $u:[N,2N]\to \r $ une fonction de classe ${\Cal C}^4$,
qui est ``petite" devant $f$. Pour cela, on suppose qu'il existe un
r\'eel $\eta \ (0\le \eta\le 1/4) $ tel que l'on ait
$$
|u^{(j)}(x)|\le \eta \lambda_2 N^{2-j}, \text{ \ pour \ }
j=1,2,3,4 \text{ \ et pour \ } x\in [N,2N]. 
\leqno (6.5)
$$
Le calcul suivant est d\'etaill\'e dans [5]:

\proclaim{Lemme 4} Soient $f$ et $u$ v\'erifiant (6.1) et (6.5). On
conserve les notations (6.3) et (6.4). 
On pose $g(x)=f(x)+u(x)$. Soit $[\alpha_0, \beta_0]$ l'intersection
des domaines de d\'efinition de $f^{*}$ et de $g^{*}$. Alors, pour
$y\in [\alpha_0,\beta_0]$, on a:
$$
g^{*}(y)=f^{*}(y)+u(z(y))-\frac{u'(z(y))^2}{2f''(z(y))}+v(y), 
$$
o\`u $v$ est une fonction de classe ${\Cal C}^1$ qui v\'erifie: \ 
$v(y)\ll \eta^3N^2\lambda_2$ et $v'(y)\ll\eta^3N$.
\hfill\hfill\qed
\endproclaim

L'exemple qui nous int\'eresse est celui de la fonction $g(x)=\alpha
x^2+\gamma x^4$, avec:
$$
0<\alpha\le \frac{1}{2}, \quad |\gamma|\le \frac{\alpha}{96N^2}
\leqno(6.6)
$$
Alors le calcul de $g^{*}$ d\'ecoule du lemme 3:
$$
g^{*}(y)=-\frac{y^2}{4\alpha}+\frac{\gamma y^4}{16 \alpha^4}
-\frac{\gamma^2y^6}{16\alpha^7}+v(y),
\leqno (6.7)
$$
avec
$$
v(y)\ll \frac{|\gamma|^3N^8}{\alpha^2}\quad \text{et}\quad
v'(y)\ll\frac{|\gamma|^3N^7}{\alpha^3}. 
$$

\head{7. Une application de la transformation B}\endhead

Le coeur de la d\'emonstration du Th\'eor\`eme 1 est le r\'esultat
suivant, analogue du lemme 5 de [1]:

\proclaim{Lemme 5} Soient $N$ et $\Delta $ deux r\'eels tels que $N\ge
16$ et $N^{-1/2}\le \Delta\le 1/4$. On pose:
$$
\align
 R(\Delta, N)&=
\int_{\Delta}^{2\Delta}\int_{-N^{-3}}^{N^{-3}} 
\bigg|\sum_{N<n\le 2N}\e(\alpha n^2+\gamma n^4)\bigg|^{6}\d\alpha\,\d\gamma,
\tag 7.1
\\ 
I(N)&=
\int_{0}^{1}\int_{-N^{-3}}^{N^{-3}} 
\bigg|\sum_{N<n\le 2N}\e(\alpha n^2+\gamma n^4)\bigg|^{6}\d\alpha\,\d\gamma,
\tag 7.2
\endalign
$$
Alors on a:
$$
R(\Delta, N)\ll \Delta (\log N)^6 \Big(I(2\Delta N)+I(4\Delta
N)\Big)
\leqno (7.3)
$$
\endproclaim

1) La preuve s'appuie sur la transformation B de Van der
Corput qu'on va appliquer \`a la somme:
$$
S(\alpha, \gamma)=\sum_{N<n\le 2N}\e(g(n)), 
$$
avec $g(n)=g_{\alpha, \gamma}(n)=\alpha n^2 +\gamma n^4$, et
$\alpha\in [\Delta, 2\Delta ], \ \gamma \in [-N^{-3},N^{-3}]$

Pour que la condition (6.6) soit satisfaite, on doit supposer que 
$$
N\ge N_0:= 96^2. 
$$
Pour d\'emontrer (7.3) dans le cas $N<N_0$, il suffit de remarquer
qu'on a toujours, sans limitation sur $N \:$
$$
I(N)\gg 1, 
\leqno (7.4)
$$
Cela provient de la deuxi\`eme in\'egalit\'e de (4.3) en remarquant
que ${\Cal B}_{2p}\big(N,\delta_1, \lambda_1 \big)$ est toujours
minor\'e par le nombre de solutions triviales du syst\`eme
correspondant.

A partir d'ici, on suppose $N\ge N_0$. Alors $g'$ est \`a valeurs dans un
intervalle fixe $J_0: = [2\Delta N-4, 8\Delta N +32 ]$. L'application
de (6.2) peut donc s'\'ecrire:
$$
S(\alpha, \gamma)\ll
\bigg| \sum_{m\in J} \frac{\e(g^{*}(m))}{g''(z(m))^{1/2}} \bigg|
+N^{1/4},
\leqno (7.5)
$$
o\`u $J=J(\alpha, \gamma) $ est un sous-intervalle de $J_0$. Mais
dans la somme du membre de droite de (7.5), on peut supprimer $O(1)$
termes sans avoir \`a modifier le terme d'erreur ; en particulier, on
peut supposer que $J$ est un sous-intervalle de $J_1: = ]2\Delta N,
8\Delta N]$, et ce, pour des raisons de commodit\'e qui
appara\^{\i}tront plus loin.

On int\`egre alors (7.5) par rapport \`a $\alpha $ et $\gamma $ et,
gr\^ace \`a une sommation d'Abel, on a finalement:
$$
R(\Delta, N) \ll \Delta^{-3}
\int_{\Delta}^{2\Delta}\int_{-N^{-3}}^{N^{-3}} 
\bigg|\sum_{m\in J}\e(g^{*}(m))\bigg|^{6}\d\alpha\,\d\gamma 
+\frac{\Delta}{N^{3/2}}.
\leqno (7.6)
$$

2) Le calcul de $g^{*}(m)$ est donn\'e par (6.7):
$$
g^{*}(m)=-\frac{m^2}{4\alpha}+\frac{\gamma m^4}{16\alpha^4}
-\frac{\gamma^2 m^6}{16\alpha^7}+v(m),
$$
o\`u $v$ est une fonction dont la variation totale est $O(1)$. En
particulier, le terme $v(m)$ peut dispara\^{\i}tre par sommation
d'Abel dans la somme (7.6). D'autre part, on pose $M=2\Delta N$.
Alors, avec la notation (3.1), et si $J$ d\'esigne n'importe quel
sous-intervalle de $J_1 \ (J_1=]2\Delta N, 8\Delta N])$, on a:
$$
\sum_{m\in J}b(m)\ll\bigg|\sum_{m\sim M}b(m)\bigg|
+\bigg|\sum_{m\sim 2M}b(m)\bigg|
$$
pour toute suite $\big(b(m)\big)_m$ de nombres complexes. On a
montr\'e que:
$$
R(\Delta, N) \ll R_1(\Delta, N)+R_2(\Delta, N) + \Delta 
\leqno(7.7)
$$
avec 
$$
R_1(\Delta, N)=
\Delta^{-3}\int_{\Delta}^{2\Delta}\int_{-N^{-3}}^{N^{-3}}
\bigg|\sum_{m\sim M}\e\left(\frac{m^2}{4\alpha}-\frac{\gamma
m^4}{16\alpha^4}
+\frac{\gamma^2 m^6}{16\alpha^7}
\right)\bigg|^6\d\alpha\,\d\gamma 
\leqno (7.8)
$$
et o\`u $R_2(\Delta, N)$ se d\'efinit de la m\^eme fa\c con, mais en
rempla\c cant $M$ par $2M$. 

On va prouver les majorations:
$$
R_1(\Delta, N)\ll \Delta (\log N)^6 I(2\Delta N) 
\text{ \ et \ } R_2(\Delta, N)\ll \Delta (\log N)^6 I(4\Delta N) 
\leqno (7.9)
$$
ce qui, compte tenu de (7.4), ach\`evera la d\'emonstration du Lemme.

3) Dans l'int\'egrale (7.8), on fait le changement de variables
$x=\dfrac{1}{4\alpha}$, $y=-\dfrac{\gamma}{16\alpha^4}$. La nouvelle
variable $(x,y)$ reste dans le rectangle
$$
\Omega=\Big[\frac{1}{8\Delta},\frac{1}{4\Delta}\Big]
\times \Big[-\frac{N^{-3}}{16\Delta^4},
\frac{N^{-3}}{16\Delta^4}\Big].
$$
On recouvre $\Omega$ par des petits rectangles de la forme
$$
\omega =[c,c+1]\times [d,d+2M^{-3}], \ c\asymp \Delta^{-1}, \
d\ll(\Delta M^3)^{-1} 
$$
Le nombre de petits rectangles n\'ecessaires est $\ll \Delta^{-2}$.
Tous calculs faits, on aboutit ainsi \`a:
$$
R_1(\Delta, N)\ll\Delta 
\max\Sb {c\asymp \Delta^{-1}}\\{ d\ll (\Delta M^3)^{-1}}\endSb
\int_{c}^{c+1}\int_{d}^{d+2M^{-3}}
\bigg|\sum_{m\sim M}\e\left(xm^2+ym^4+\frac{4y^2}{x}m^6
\right)\bigg|^6\d x\d y 
$$
Maintenant, on pose 
$$
a_m=\e\Bigl(4\frac{d^2m^6}{c}\Bigr) \text{ \ et \ }
\phi_{x,y}(m)=4\Big(\frac{y^2}{x}-\frac{d^2}{c}\Big)m^6.
$$
On obtient:
$$
R_1(\Delta, N)\ll\Delta \max_{c,d}
\int_{c}^{c+1}\int_{d}^{d+2M^{-3}}
\bigg|\sum_{m\sim M}a_m \e\big(xm^2+ym^4+\phi_{x,y}(m)
\big)\bigg|^6\d x\d y 
$$
La variation totale de $\phi_{x,y}$ est $\ll 1 ;$ les hypoth\`eses du
lemme 1 sont v\'erifi\'ees et on en d\'eduit la premi\`ere
in\'egalit\'e de (7.9). Raisonnant de m\^eme pour la deuxi\`eme, on
ach\`eve la d\'emonstration du Lemme. \hfill\hfill\qed

\subhead{D\'emonstration du Th\'eor\`eme 1}\endsubhead
Nous commen\c cons par \'etablir une relation de r\'ecurrence:
\proclaim{Lemme 6} Pour $N\ge 16$, on d\'efinit $I(N)$ comme en
(7.2). On suppose qu'il existe un r\'eel $\beta\ge 0$ tel qu'on ait:
$$
I(N)\ll_{\ep}N^{\beta+\ep}
\leqno (7.10)
$$
Alors on a:
$$
I(N)\ll_{\ep}N^{\frac{\beta}{1+\beta}+\ep}
\leqno (7.11).
$$
\endproclaim

On pose $\Delta=N^{-\frac{\beta}{1+\beta}}$. Par le Lemme 1, on a:
$$
I(N)\ll \Delta^{-1}(\log N)^6\int_0^{\Delta}\int_{-N^{-3}}^{N^{-3}}
\bigg|\sum_{n}\e(\alpha n^2 + \gamma n^4)
\bigg|^{6}\d\alpha\,\d\gamma
$$
On d\'ecompose l'int\'egrale: $\int_0^\delta \le \int_0^{N^{-1/2}}+
\sum_{N^{-1/2}\le \delta <\Delta}\int_{\delta}^{2\delta}$, o\`u
$\delta$ est de la forme $2^kN^{-1/2}$. On utilise le Lemme 3 pour
l'int\'egrale $\int_0^{N^{-1/2}}$, d'o\`u
$$
I(N)\ll \Delta^{-1}(\log N)^{12}+\Delta^{-1}(\log
N)^{7}\max_{N^{-1/2}\le \delta\le\Delta}\big|R(\delta,N)\big|,
$$
o\`u $R(\delta,N)$ est d\'efini en (7.1). On applique le Lemme 5, puis
l'hypoth\`ese de r\'ecurrence (7.10) et on obtient:
$$
I(N)\ll_{\ep}\Delta^{-1}N^{\ep}+\Delta^{\beta}N^{\beta+\ep},
$$
ce qui, par le choix de $\Delta, $ redonne (7.11).\hfill\hfill\qed

{\it D\'emonstration du Th\'eor\`eme 1:}
 On doit majorer 
$$
I_6 =\int_0^{1}\int_{0}^{\lambda}
\bigg|\sum_{n\sim N}\e(\alpha n^2 + \gamma n^4)
\bigg|^{6}\d\alpha\,\d\gamma. 
$$
En consid\'erant s\'epar\'ement le cas $\lambda\le N^{-3}$ et le cas
$\lambda>N^{-3}$, le Lemme 1 montre qu'on a:
$$
I_6\ll(\log N)^6(1+\lambda N^3)I(N). 
$$
D'autre part, une majoration triviale montre que (7.10) est vrai pour
$\beta=\beta_0:= 3$. On d\'efinit par r\'ecurrence, pour $n\ge 1$,
$\beta_n=\frac{\beta_{n-1}}{1+\beta_{n-1}}$. Par le lemme 6 appliqu\'e
$n$ fois, on obtient (7.10) avec $\beta=\beta_n$, ce qui donne, en
prenant $n$ arbitrairement grand:
$$
I(N)\ll_{\ep}N^{\ep},
$$
et le th\'eor\`eme 1 en d\'ecoule. \hfill\hfill\qed

Pour en finir avec les puissances sixi\`emes, nous rappelons la
formulation \'equiva\-lente du Th\'eor\`eme 1 en termes de syst\`emes
diophantiens:

\proclaim{Th\'eor\`eme 1'} Soit un entier $N\ge 2$. Le nombre de
solutions du syst\`eme
$$
\aligned
&n_1^2+n_2^2+n_3^2 = m_1^2+m_2^2+m_3^2\\
&|n_1^4+n_2^4+n_3^4 - m_1^4-m_2^4-m_3^4| \le \delta N^4
\endaligned\quad
0\le n_i,m_i\le N
$$
est \  
$
\ll_{\ep} N^{3+\ep}+\delta N^{4+\ep}.
$
\hfill\hfill\qed
\endproclaim
\head{8. D\'emonstration du Th\'eor\`eme 2}\endhead

On doit majorer
$$
I_8=\int_0^1\int_0^{\lambda}\bigg|\sum_{n\sim N}a_n\e(\alpha
n^2+\gamma n^4)\bigg|^8\d\alpha\,\d\gamma.
$$
On pose $\mu=N^{-5/2}$. En s\'eparant le cas $\lambda\le\mu$ et le cas
$\lambda>\mu$, on obtient, par le Lemme 1:
$$
I_8\ll_\ep N^\ep \Big(1+\frac{\lambda}{\mu}\Big)
\int_0^1\int_0^{\mu}\bigg|\sum_{N<n\le 2N}\e(\alpha
n^2+\gamma n^4)\bigg|^8\d\alpha\,\d\gamma.
$$

On d\'ecompose l'int\'egrale: 
$$
\int_0^\mu\le\int_0^{N^{-3}}+\sum_{N^{-3}\le \delta < \mu}
\int_\delta^{2\delta},
$$
avec $\delta$ de la forme: $\delta = 2^k N^{-3}$, d'o\`u:
$$
I_8\ll_\ep N^\ep \Big(1+\frac{\lambda}{\mu}\Big)
\Big(J_0+\max_{N^{-3}\le \delta \le \mu}J_\delta \Big),
\leqno (8.1)
$$
o\`u on a pos\'e:
$$
\align
J_0&=\int_0^1\int_0^{N^{-3}}\bigg|\sum_{n}\e(\alpha
n^2+\gamma n^4)\bigg|^8\d\alpha\,\d\gamma 
\\
J_\delta &=\int_0^1\int_\delta ^{2\delta}\bigg|\sum_{n}\e(\alpha
n^2+\gamma n^4)\bigg|^8\d\alpha\,\d\gamma.
\endalign
$$
La premi\`ere int\'egrale se majore simplement:
$$
J_0\ll N^2\int_0^1\int_0^{\lambda}\bigg|\sum_{n}\e(\alpha
n^2+\gamma n^4)\bigg|^6\d\alpha\,\d\gamma\ll_\ep N^{2+\ep},
$$
d'apr\`es le Th\'eor\`eme 1. Pour les autres int\'egrales, on \'ecrit
:
$$
\align
J_\delta &\le \bigg(
\max\Sb 0\le\alpha\le 1\\\delta\le\gamma\le 2\delta\endSb
\bigg|\sum_{N<n\le 2N}\e(\alpha n^2+\gamma n^4)\bigg|\bigg)^2
\int_0^1\int_0^{2\delta}\Big|\sum \Big|^6 \d\alpha\,\d\gamma 
\\
&\ll_\ep \Bigg(N^2(\delta N)^{1/3}+N^{3/2}+\delta^{-1/2}\Bigg)\delta
N^{3+\ep},
\endalign
$$
d'une part d'apr\`es le crit\`ere de la d\'eriv\'ee troisi\`eme pour
les sommes d'exponentielles (cf\. [2, Th\'eor\`eme 2.6]) et, d'autre
part, d'apr\`es le Th\'eor\`eme 1. Par le choix de $\mu$, on a
$J_\delta \ll_\ep N^{2+\ep}$ pour tout $\delta\in [N^{-3},\mu]$.
Reportant ces majorations dans (8.1), on obtient (3.2). \hfill\hfill\qed

\head{9. D\'emonstration du Th\'eor\`eme 3}\endhead
 
Le Th\'eor\`eme 3 se d\'eduit du Th\'eor\`eme 2 exactement de la
m\^eme fa\c con que le Th\'eor\`eme 2 se d\'eduit du Th\'eor\`eme
1. La d\'emonstration pr\'ec\'edente peut \^etre reconduite avec les
modifications suivantes:

Tout d'abord, on doit prendre $\mu=N^{-17/8}$ au lieu de $N^{-5/2}$.
Le deuxi\`eme point \`a changer est le d\'ecoupage de l'int\'egrale
qu'on doit maintenant \'ecrire:
$$
\int_0^\mu \le \int_0^{N^{-5/2}}+\sum_{N^{-5/2}\le \delta \le
\mu}\int_{\delta}^{2\delta },
$$
toujours avec $\delta$ de la forme $\delta=2^kN^{-5/2}$.

Enfin, chaque application du Th\'eor\`eme 1 doit \^etre remplac\'ee
par une application du Th\'eo\-r\`eme~2. \hfill\hfill\qed

{\it Remarques:}
1) Dans la d\'emonstration du Th\'eor\`eme 3, on doit majorer la somme
d'exponentielles:
$$
 S(\mu)=\max\Sb 0\le \alpha\le 1\\ \gamma \asymp \mu\endSb
\Big|\sum_n\e(\alpha n^2+\gamma n^4)\Big|, \text{ \ avec \ } \mu=N^{-17/8}.
$$
 Le crit\`ere de la d\'eriv\'ee troisi\`eme
s'applique pr\'ecis\'ement \`a ce type de situation. Cependant, on
peut faire mieux en introduisant des techniques plus \'elabor\'ees(par
exemple, en appliquant le Lemme A, puis en adaptant une paire
d'exposants convenable, cf\. [2]). Mais les am\'eliorations attendues
sont minimes ; c'est pourquoi nous n'avons pas cherch\'e \`a le faire.

2) La m\^eme m\'ethode permet de d\'eduire du Th\'eor\`eme 3 un
r\'esultat sur les puissances douzi\`emes. Il faudrait pour cela
prendre $\mu$ nettement plus grand que $N^{-17/8}$. Dans ce cas, les
majorations possibles de $S(\mu)$, m\^eme \`a l'aide de paires
d'exposants performantes, fourniraient un r\'esultat insuffisant pour
l'application aux sommes de Weyl.

\head{10. Application \`a l'in\'egalit\'e de Weyl}\endhead

{\bf 10.1. Enonc\'e du r\'esultat.}
Soit $k$ un entier fix\'e, $k\ge 8 ;$ les constantes sous-entendues
dans le symbole $\ll$ pourront dor\'enavant d\'ependre de $k$.
Etant donn\'e un r\'eel $\alpha$ et un entier $N\ge 2$, on veut
majorer la somme d'exponentielles
$$
S(\alpha)=S_{N,k}(\alpha)=\sum_{n=1}^N\e(\alpha n^k) 
$$
\`a l'aide de la quantit\'e
$
B_{H,\delta}(\alpha)=\# \big\{h=1,\dots,H \;\big|\; \|\alpha h\|\le \delta
\big\}, 
$
o\`u $H$ est un entier, o\`u $\delta$ est un r\'eel positif et o\`u
$\|x\|$ d\'esigne la distance du r\'eel $x$ \`a l'entier le plus
proche.

\proclaim{Th\'eor\`eme 4} On pose $H=16\frac{k!}{4!} N^{k-4}, \ 
\delta=N^{-4}$ et $K=2^k. $ Alors on a:
$$
S(\alpha)\ll_\ep N^{1-\frac{16}{K}}+N^{1-\frac{3}{K}+\ep}
\Big(\frac{B_{H,\delta}(\alpha)}{HN^{-4}}
\Big)^{\frac{8}{5K}}. 
\leqno (10.1)
$$
\endproclaim\medskip
{\bf 10.2. Remarques sur le r\'esultat.}
Supposons que $B_{H,\delta}(\alpha)$ admette la majoration
``probabiliste":
$$
B_{H,\delta}(\alpha)\ll_\ep H^{1+\ep}\delta +H^{\ep}. 
\leqno
(10.2)
$$
Alors (10.1) devient:
$$
S(\alpha)\ll_{\ep} N^{1+\ep-\frac{3}{K}}\big(1+N^{8-k}\big)^{8/5K},
$$
ce qui explique pourquoi nous nous restreignons au cas $k\ge 8$. 
Par ailleurs, on peut majorer $B_{H,\delta}(\alpha)$ en fonction des
approximations rationnelles de $\alpha$ gr\^ace \`a un lemme de
Heath-Brown que nous rappelons maintenant:

\proclaim{Lemme 7} (cf\. [3, Lemme 6]). On suppose que $\alpha $ admet
l'approximation rationnelle
$$
\alpha =a/q+\theta, \text{ \ avec \ } q\ge 1, \ (a,q)=1 
\text{ \ et \ } |\theta |\le q^{-2}. 
\leqno (10.3)
$$
Alors on a les deux majorations:
$$
 B_{H,\delta}(\alpha)\le 4(1+q\delta\big)\big(1+H/q)
\quad
\text{et} 
\quad
 B_{H,\delta}(\alpha)\le 8\big(1+{\delta}/{q|\theta
 |}\big)\big(1+q|\theta|H\big). 
$$
\hfill\hfill\qed
\endproclaim
On en d\'eduit le:
\proclaim{Corollaire} Avec les notations (10.3), on suppose que l'une
des deux hypoth\`eses suivantes est v\'erifi\'ee:
$$
N^4\ll q\ll N^{k-4}
\quad
\text{ou}
\quad
\frac{1}{N^{k-4}}\ll |\theta|q \ll \frac{1}{N^4}. 
$$
Alors, on a la majoration \ 
$
S(\alpha)\ll_{\ep} N^{1+\ep-\frac{3}{K}}
$
\endproclaim\medskip
{\bf 10.3. Deux lemmes classiques.}
Le premier lemme est le Lemme A it\'er\'e $r$ fois faisant intervenir
les diff\'erences sym\'etriques: $\Delta_hf(x)=f(x+h)-f(x-h)$. Pour
une d\'emonstration, nous renvoyons au $\S 2$ de [3].
\proclaim{Lemme 8} Pour toute fonction $f:[1,N]\to \r$, et tout
entier $r\ge 1$, on a, en posant $R=2^r:$
$$
\multline
 \Big|\sum_{n=1}^N\e(f(n))\Big|^R\ll_r\\
 N^{R-1} 
+N^{R-(r+1)}\sum_{i=1,2}\sum\Sb|h_1|<N/2\\h_1\not = 0\endSb
\dots\sum\Sb|h_r|<N/2\\h_r\not = 0\endSb
\bigg|\sum\Sb n\in J(\bs h)\\ n\equiv i\pmod2\endSb
 \e(\Delta_{h_1}\dots\Delta_{h_r}f(n))\bigg|, 
\endmultline
$$
o\`u $J({\bs h})=J(h_1,\dots h_r)$ est un sous-intervalle de
$[1,N]$. \hfill\hfill\qed
\endproclaim
Le deuxi\`eme lemme est un cas particulier du double grand crible de
Bombieri et Iwaniec (cf\. [2, Lemme 7.5]):

\proclaim{Lemme 9} Pour tout entier $h=1,\dots,H$, soient $a_h$ et
$b_h$ des r\'eels. Soit un entier fix\'e $p\ge 1$. Pour chaque $h$,
soit $N_h$ un entier tel que $1\le N_h\le N$. Alors on a:
$$
\multline
\sum_{h=1}^H\bigg|\sum_{n=1}^{N_h}\e (a_h n^4+b_h n^2)\bigg|\\
\ll_{\ep} H^{1-\frac{1}{p}} N^{3/p +\ep}\Cal N^{1/2p}
\bigg(\int_0^1\int_0^1\Big|\sum_{n=1}^N
\e(xn^2+yn^4)\Big|^{2p}\d x\d y\bigg)^{1/2p} 
\endmultline
$$
avec
$\Cal N\!=\!\#\big\{ (h_1,h_2)\in \{1,\dots,H\}^2 \,\big|\,
\|a_{h_1}\!-\!a_{h_2}\|\le
N^{-4} \text{ et } \|b_{h_1}\!-\!b_{h_2}\|\le N^{-2} \big\}$.\qed
\endproclaim
\medskip
{\bf 10.4. D\'emonstration du th\'eor\`eme 4.}
On applique le Lemme A it\'er\'e $k-4$ fois (au lieu de $k-3$ fois
dans [3]). Par le calcul, on obtient:
$$
\Delta_{h_1}\dots\Delta_{h_r}(\alpha n^k)= 2^{k-4}\frac{k!}{4!} h_1
\dots h_{k-4}\alpha n^4+b_{\underline {h}}n^2+c_{{\bs h}}, 
$$
l'expression de $b_{{\bs h}}=b_{h_1,\dots,h_{k-4}}$ et de
 $c_{{\bs h}}=c_{h_1,\dots,h_{k-4}}$ \'etant sans importance. 

On doit maintenant r\'ealiser quelques op\'erations triviales pour se
ramener exactement au Lemme 9. 

On pose $h=2^{k-4}\frac{k!}{4!}h_1\dots h_{k-4}$. On d\'esigne ensuite par
$b_h'$ le r\'eel qui r\'ealise le maximum de la quantit\'e:
$$
\max_{1\le N_1 \le N}\bigg|\sum _{n=1}^{N_1}\e(h\alpha n^4+b_h'
n^2)\bigg|.
$$
Enfin, on remarque que, pour tout intervalle $J\subset [1,N]$, on a:
$$
\sum_{i=1,2}\Big|\sum_{n\in J}\beta(n)\Big|\ll
\max_{1\le N_1\le N}\Big|\sum_{n=1}^{N_1}\beta(n)\Big|+
\max_{1\le N_1\le N}\Big|\sum_{n=1}^{N_1/2}\beta(2n)\Big|.
$$
pour toute suite $(\beta(n))_n$ de nombres complexes.
\smallskip
On a finalement montr\'e que, pour un bon choix des entiers $N_h$, on
a:
$$
\big|S(\alpha)\big|^{K/16}\ll_{\ep} N^{K/16-1}
+ N^{K/16-k+3+\ep}\sum_{h=1}^H\Big|\sum_{n=1}^{N_h}
\e(h\alpha n^4 +b_h' n^2)\Big|, 
$$
o\`u $H$ est d\'efini comme au Th\'eor\`eme 4.
\smallskip
On applique maintenant le Lemme 9 avec $p\!=\!5$, puis on majore 
$\Cal N$ par

 $B_{H,N^{-4}}(\alpha),$ ce qui donne (10.1) gr\^ace au
Th\'eor\`eme 3. \hfill\hfill\qed
\smallskip{\it Remarque: } Si on utilise le Th\'eor\`eme 2 \`a la place du
Th\'eor\`eme 3, on obtient
$$
S(\alpha)\ll N^{1-16/K}+N^{1-3/K+\ep}
\Big(\frac{B_{H,\delta}(\alpha)}{HN^{-4}}\Big)^{2/K}. 
$$
Ce r\'esultat est identique \`a (10.1) dans le cas d'une majoration
probabiliste de $B_{H,\delta}(\alpha)$ (cf (10.2)), mais il est moins
bon dans les autres cas.

\Refs

\item{[1]} E. Bombieri et H. Iwaniec, {\it Some mean value theorems for
exponential sums, }Ann. Scuola Norm. Sup. Pisa Cl. Sci (4) {\bf 13}
(1986), 473--486
\item{[2]} S. W. Graham and G. Kolesnik,  {\it Van der Corput's Method
for Exponential Sums, } London Math. Soc. Lecture Notes Series 126,
Cambridge University Press, 1991

\item{[3]} D. R. Heath-Brown,  {\it Weyl's inequality, Hua's
inequality, and Waring's problem, } J. London Math. Soc. (2) {\bf 38}
(1988), 216--230

\item{[4]} H. L. Montgomery, {\it Ten Lectures on the Interface Between
Analytic Number Theory and Harmonic Analysis, }Conference board of the
Mathematical Sciences,  American Mathematical, Society, 1994

\item{[5]} M. Redouaby et P. Sargos,  {\it Sur la transformation B de
Van der Corput, } Expo. Math. {\bf 17} (1999), 207--232

\item{[6]} P. Sargos, {\it Un crit\`ere de la d\'e\-riv\'ee
cin\-qui\`eme pour les som\-mes d'expo\-nen\-tielles, }\`a para\^\i tre
\`a Bull. London Math. Soc.

\item{[7]} R. C. Vaughan, {\it The Hardy-Littlewood Method, } Cambridge
University Press, 1981

\item{[8]} N. Watt,  {\it Exponential sums and the Riemann zeta
function II, } J. London Math. Soc. (2) {\bf 39} (1989), 385--404

\endRefs

\enddocument